\documentclass[11pt]{amsart}
\usepackage{CJK}
\usepackage{amssymb}
\usepackage{amsmath}
\usepackage{amsthm}
\usepackage{bm}

\usepackage{graphicx}

\newtheorem{Prop}{Proposition}[section]
\newtheorem{Thm}[Prop]{Theorem}
\newtheorem{Lem}[Prop]{Lemma}
\newtheorem{Cor}[Prop]{Corollary}

\newtheorem*{Thm0}{Theorem}

\theoremstyle{definition}
\newtheorem{Ex}[Prop]{Example}

\begin{document}

\title{A FURTHER NOTE ON THE CONCORDANCE INVARIANTS EPSILON AND UPSILON}

%\date{March 2016}
\author{Shida Wang}

\address{Department of Mathematics, University of Oregon, Eugene, OR 97403}
\email{wang217@indiana.edu}

\begin{abstract}

Hom gives an example of a knot with vanishing Upsilon invariant but nonzero epsilon invariant.
We build more such knots that are linearly independent in the smooth concordance group.

\end{abstract}

\maketitle

\section{Introduction}

Let $\mathcal{C}$ be the smooth concordance group. It has a subgroup $\mathcal{C}_{TS}$ consisting of topologically slice knots.
The $\varepsilon$ invariant~\cite{epsilon} and $\Upsilon$ invariant~\cite{upsilon} derived from knot Heegaard Floer theory have shown their power in proving the following result.

\begin{Thm0}\emph{(\cite[Theorem 1]{1nn1} and \cite[Theorem 1.20]{upsilon})} The group $\mathcal{C}_{TS}$ contains a summand isomorphic to $\mathbb{Z}^\infty$.\end{Thm0}

One may wonder whether one of the two invariants is stronger than the other. In fact, Hom shows that the $\varepsilon$ invariant is not weaker than the $\Upsilon$ invariant.

\begin{Thm0}\emph{(\cite[Theorem 2]{original})} There exists a knot with vanishing $\Upsilon$ invariant but nonvanishing $\varepsilon$ invariant.\end{Thm0}

It is not known yet whether there is a knot with vanishing $\varepsilon$ invariant but nonvanishing $\Upsilon$ invariant.

We will prove the following result.

\begin{Thm}\label{main}There exists a subgroup of $\mathcal{C}$ isomorphic to $\mathbb{Z}^\infty$ such that each of its nonzero elements has vanishing $\Upsilon$ invariant but nonvanishing $\varepsilon$ invariant.\end{Thm}

Moreover, we will give many such subgroups.

The knots generating the subgroup will be built from known examples that have vanishing $\Upsilon$ invariant.
For any knot $K$, the $\Upsilon$ invariant, denoted by $\Upsilon_K(t)$, is a piecewise linear function on $[0,2]$.
%We refer to \cite{upsilon} for the definition and basic properties of $\Upsilon$ invariant.
This invariant gives a homomorphism from $\mathcal{C}$ to the vector space of continuous functions on $[0,2]$.
Feller and Krcatovitch show the following recursive formula for the $\Upsilon$ invariant of torus knots.

\begin{Thm}\label{expansion}\emph{(\cite[Proposition 2.2]{recursive})} Suppose $p$ and $q$ are relatively prime positive integers and $k$ is a nonnegative integer.
Then $\Upsilon_{T_{q,kq+p}}(t)=\Upsilon_{T_{p,q}}(t)+k\Upsilon_{T_{q,q+1}}(t)$.\end{Thm}

From this formula, one immediately knows that any knot of the form $T_{q,kq+p}\#-T_{p,q}\#-kT_{q,q+1}$ has vanishing $\Upsilon$ invariant,
where $-K$ means the mirror image of the knot $K$ with reversed orientation (representing the inverse element of $K$ in $\mathcal{C}$) and $kK$ means the connected sum of $k$ copies of $K$.
By a theorem of Litherland~\cite{signature} that torus knots are linearly independent in the concordance group,
it is easy to give subgroups of $\mathcal{C}$ isomorphic to $\mathbb{Z}^\infty$ included in the kernel of the $\Upsilon$ homomorphism,
but its elements may also have vanishing $\varepsilon$ invariant~\cite[Remark 4.11]{secondary2}.
In Proposition~\ref{three torus knots}, we will show a bound for some of the knots of the form $T_{q,kq+p}\#-T_{p,q}\#-kT_{q,q+1}$ with respect to the total order given by the $\varepsilon$ invariant,
which enables us to prove Theorem~\ref{main}.
Such a bound also yields examples with arbitrarily large concordance genus but vanishing $\Upsilon$ invariant.

\emph{Acknowledgments.} The author wishes to express sincere thanks to Professor Robert Lipshitz for carefully reading a draft of this paper.

\section{Preliminaries}

\subsection{The $\varepsilon$ invariant}

We assume the reader is familiar with knot Floer homology, defined by Ozsv\'{a}th-Szab\'{o}~\cite{CFKdef1} and independently Rasmussen~\cite{CFKdef2},
and the $\varepsilon$ invariant, defined by Hom~\cite{epsilon}.
We briefly recall some properties of the $\varepsilon$ invariant and its refinement for later use.

To a knot $K\subset S^3$, the knot Floer complex associates a doubly filtered, free, finitely generated chain complex over $\mathbb{F}[U,U^{-1}]$, denoted by $\mathit{CFK}^\infty(K)$,
where $\mathbb{F}$ is the field with two elements.
Up to filtered chain homotopy equivalence, this complex is an invariant of $K$.
Hom defines an invariant of the filtered chain homotopy type of $\mathit{CFK}^\infty(K)$ and hence an invariant of $K$, called $\varepsilon$, taking on values $-1$, 0 or~1~\cite{epsilon},
which has the following properties~\cite[Proposition 3.6]{epsilon}:

(1) if $K$ is smoothly slice, then $\varepsilon(K)=0$;\\
(2) $\varepsilon(-K)=-\varepsilon(K)$;\\
(3) if $\varepsilon(K)=\varepsilon(K')$, then $\varepsilon(K\#K')=\varepsilon(K)=\varepsilon(K')$;\\
(4) if $\varepsilon(K)=0$, then $\varepsilon(K\#K')=\varepsilon(K')$.

Thus the relation $\sim$, defined by $K\sim K'\Leftrightarrow\varepsilon(K\#-K')=0$, is an equivalence relation coarser than smooth concordance.
It gives an equivalence relation on $\mathcal{C}$ called \emph{$\varepsilon$-equivalence}.
The~$\varepsilon$-equivalence class of $K$ is denoted by $[\![K]\!]$.
All $\varepsilon$-equivalence classes form an abelian group, denoted by~$\mathcal{CFK}$ in~\cite{1nn1}, which is a quotient group of $\mathcal{C}$.
It follows from the properties of~$\varepsilon$ that the definition $[\![K]\!]>[\![K']\!]\Leftrightarrow\varepsilon(K\#-K')=1$
gives a total order on $\mathcal{CFK}$ that respects the addition operation~\cite[Proposition 4.1]{ordered}.

More generally, the $\varepsilon$ invariant can be defined for a larger class $\mathfrak{C}$ of doubly filtered, free, finitely generated chain complex over $\mathbb{F}[U,U^{-1}]$ with certain homological conditions
(see~\cite[Definition 2.2]{1nn1}),
which contains $\mathit{CFK}^\infty(K)$ for all knots $K$.
Similar properties hold, with the connected sum operation replaced by tensor product operation, and the negative of a knot corresponds to the dual of a complex.
The~$\varepsilon$-equivalence class of $C$ in $\mathfrak{C}$ is denoted by $[C]$.
All $\varepsilon$-equivalence classes of complexes in $\mathfrak{C}$ also form a totally ordered abelian group, denoted by~$\mathcal{CFK}_\mathrm{alg}$ in~\cite{1nn1},
which includes~$\mathcal{CFK}$ as a subgroup.

For any complex $C$ in $\mathfrak{C}$ with $\varepsilon(C)=1$, Hom defines a tuple of numerical invariants $\bm{a}(C)=(a_1(C),\cdots,a_n(C))$ in~\cite[Section 3]{1nn1}, where the positive integer $n$ depends on $C$.
For a knot $K$ with $\varepsilon(K)=1$, denote $\bm{a}(\mathit{CFK}^\infty(K))=\bm{a}(K)=(a_1(K),\cdots,a_n(K))$.
The numbers $a_1(C),\cdots,a_n(C)$ for any complex $C$ in $\mathfrak{C}$ satisfy exactly one of the following three conditions:

(1) $a_1(C),\cdots,a_n(C)$ are positive integers;\\
(2) $a_1(C),\cdots,a_{n-1}(C)$ are positive integers, $n>1$ and $a_n(C)$ is a negative integer less than~$-1$;\\
(3) $a_1(C),\cdots,a_{n-2}(C)$ are positive integers, $n>2$, $a_{n-1}(C)=-1$ and $a_n(C)$ is a negative integer.

\subsection{The staircase complex and semigroups of torus knots}

There is a special type of complexes in $\mathfrak{C}$, called \emph{staircase complexes}, that can be described by lengths of differential arrows (see~\cite[Section 2.4]{filtration} for example).
A staircase complex is usually encoded by an even number of positive integers $b_1,\cdots,b_{2m}$ that are palindromic, meaning $b_i=b_{2m+1-i}$ for $i=1,\cdots,2m$.
We will denote such a complex by $\mathrm{St}(b_1,\cdots,b_{2m})$.
Using the definition of $\bm{a}$, it can be verified that $\bm{a}(\mathrm{St}(b_1,\cdots,b_{2m}))=(b_1,\cdots,b_{2m})$.
We denote the $\varepsilon$-equivalence class of the complex $\mathrm{St}(b_1,\cdots,b_{2m})$ by $[b_1,\cdots,b_{2m}]$.
Note that this notation coincides that on~\cite[pp. 1087--1088]{1nn1}, but it is denoted by just $[b_1,\cdots,b_m]$ in~\cite{filtration}.

An example is $\mathit{CFK}^\infty(K)$ for any $L$-space knot $K$, where $b_1,\cdots,b_{2m}$ are determined by the Alexander polynomial~$\Delta_K(t)$~\cite[Remark 6.6]{ordered}.
Here $\Delta_K(t)$ must be of the form $\sum_{i=0}^{2m}(-1)^it^{\alpha_i}$, where $0=\alpha_0<\alpha_1<\cdots<\alpha_{2m}$~\cite[Theorem 1.2]{Lknot}.

Any torus knot $T_{p,q}$ is an $L$-space knot. To see which staircase complex $\mathit{CFK}^\infty(T_{p,q})$ is in terms of $p$ and $q$, it is convenient to use the notion of semigroups.
Denote the semigroup $\langle p,q\rangle=\{px+qy\mid x,y\in\mathbb{Z}_{\geqslant0}\}$ by $S$ and suppose that $\mathit{CFK}^\infty(T_{p,q})=\mathrm{St}(b_1,\cdots,b_{2m})$.
Then $b_1,\cdots,b_{2m}$ are determined by $S$ in the following way:
\begin{align*}
0,\cdots,b_1&\in S,\\
b_1+1,\cdots,b_1+b_2&\not\in S,\\
b_1+b_2+1,\cdots,b_1+b_2+b_3&\in S,\\
b_1+b_2+b_3+1,\cdots,b_1+b_2+b_3+b_4&\not\in S,\\
&\ \vdots\\
n&\in S, \forall n\geqslant b_1+\cdots+b_{2m}.
\end{align*}
This is because $S$ determines $\Delta_{T_{p,q}}(t)$ by the equation $(1-t)(\sum_{s\in S}t^s)=\Delta_{T_{p,q}}(t)$~\cite{Puiseux1} and~$\Delta_{T_{p,q}}(t)$ determines $b_1,\cdots,b_{2m}$.

Finally, the $\varepsilon$-equivalence classes of staircase complexes can be decomposed in $\mathcal{CFK}_\mathrm{alg}$ sometimes.

\begin{Lem}\label{split}\emph{(\cite[Lemma 3.1]{filtration})} Let $a_i,b_j$ be positive integers for $i=1,\cdots,k$ and $j=1,\cdots,m$.
If $k$ is even and $\max\{a_i\mid i\text{ is odd}\}\leqslant b_j\leqslant\min\{a_i\mid i\text{ is even}\}$ for all $j=1,\cdots,m$,
then $$[a_1,\cdots,a_k,a_k,\cdots,a_1]+[b_1,\cdots,b_m,b_m,\cdots,b_1]=[a_1,\cdots,a_k,b_1,\cdots,b_m,b_m,\cdots,b_1,a_k,\cdots,a_1].$$\end{Lem}

\begin{Ex}Abbreviate the finite sequence $\underbrace{1,n,\cdots,1,n}_k$ to $(1,n)^k$.
If the positive integer $b_j\leqslant n$ for every $j=1,\cdots,m$, then $$[1,n,b_1,\cdots,b_m,b_m,\cdots,b_1,n,1]=[1,n,n,1]+[b_1,\cdots,b_m,b_m,\cdots,b_1].$$
By induction, $[(1,n)^k,b_1,\cdots,b_m,b_m,\cdots,b_1,(n,1)^k]=k[1,n,n,1]+[b_1,\cdots,b_m,b_m,\cdots,b_1]$.
In particular, $[(1,n)^k,(n,1)^k]=k[1,n,n,1]$.\end{Ex}

\subsection{The totally ordered abelian group}

Let $G$ be a totally ordered abelian group, that is, an abelian group with a total order respecting the addition operation. Denote its identity element by 0.
For two elements $g,h\geqslant0$ of $G$, we write $g\ll h$ if $N\cdot g<h$ for any natural number $N$.
The \emph{absolute value} of an element $g\in G$ is defined to be $|g|=\left\{\begin{aligned}&g&\text{ if }g\geqslant0,\\&-g&\text{ if }g<0.\end{aligned}\right.$
We say that $h$ \emph{dominates} $g$ if $|g|\ll|h|$.
%Any nonzero element is defined to dominate $0$.

\begin{Lem}\label{domination implies independence}\emph{(\cite[Lemma 4.7]{ordered})} If $0<g_1\ll g_2\ll g_3\ll\cdots$ in $G$, then $g_1,g_2,g_3,\cdots$ are linearly independent in $G$.\end{Lem}

The invariants $a_1$ and $a_2$ are useful in determining domination.

\begin{Lem}\label{a1a2greater}\emph{(\cite[Lemmas 6.3 and 6.4]{ordered})} If $\bm{a}(C)=(a_1(C),\cdots)$ and $\bm{a}(C')=(a_1(C'),\cdots)$ with $a_1(C)>a_1(C')>0$,
then $[C]\ll[C']$.

If $\bm{a}(C)=(a_1(C),a_2(C),\cdots)$ and $\bm{a}(C')=(a_1(C'),a_2(C'),\cdots)$ with $a_1(C)=a_1(C')>0$ and $a_2(C)>a_2(C')>0$,
then $[C]\gg[C']$.\end{Lem}

\section{Proof of the result}

\begin{Lem}\label{remainder1} For any torus knot $T_{p,q}$ with $4\leqslant p<q$, if $q=kp+1$ for a positive integer $k$, then $[\![T_{p,q}]\!]=k[1,p-1,p-1,1]+O$
where $O\in\mathcal{CFK}_\mathrm{alg}$ is dominated by $[1,n,n,1]$ for any positive integer $n$.\end{Lem}

\textbf{Proof.} The semigroup $\langle p,q\rangle=\{0,p,2p,\cdots,kp,kp+1,kp+p,\cdots\}$.
This implies that $\mathit{CFK}^\infty(T_{p,q})=\mathrm{St}((1,p-1)^k,2,\cdots)$.
It follows from~\cite[Lemma 4.8]{1nn1} that there is some $C\in\mathfrak{C}$ with $\bm{a}(C)=(2,\cdots)$ such that $[\![T_{p,q}]\!]=[(1,p-1)^k,(p-1,1)^k]+[C]$.
Then Lemma~\ref{a1a2greater} gives the conclusion.
\hfill$\Box$

\begin{Lem}\label{any remainder} For any torus knot $T_{p,q}$ with $4\leqslant p<q$, if $q=kp+r$ for positive integers $k$ and $r$ such that $r<p$, then $[\![T_{p,q}]\!]=k[1,p-1,p-1,1]+O$
where $O\in\mathcal{CFK}_\mathrm{alg}$ is dominated by $[1,p-1,p-1,1]$.\end{Lem}

\textbf{Proof.} If $r=1$, this directly follows from Lemma~\ref{remainder1}.

Now suppose $r>1$. Then the semigroup $\langle p,q\rangle=\{0,p,2p,\cdots,kp,kp+r,kp+p,\cdots\}$.
This implies that $\mathit{CFK}^\infty(T_{p,q})=\mathrm{St}((1,p-1)^k,1,r-1,\cdots)$.
It follows from~\cite[Lemma 4.5]{1nn1} that there is some $C\in\mathfrak{C}$ with $\bm{a}(C)=(1,r-1,\cdots)$ such that $[\![T_{p,q}]\!]=[(1,p-1)^k,(p-1,1)^k]+[C]$.
Then Lemma~\ref{a1a2greater} gives the conclusion.
\hfill$\Box$

\begin{Lem}\label{half remainder} For any torus knot $T_{p,q}$ with $p<q$, if $q=kp+r$ for positive integers $k$ and $r$ such that $3\leqslant r<p/2$,
then $[\![T_{p,q}]\!]=k[1,p-1,p-1,1]+O$ where $O\gg[1,r-1,r-1,1]$ in $\mathcal{CFK}_\mathrm{alg}$.\end{Lem}

\textbf{Proof.} The semigroup $\langle p,q\rangle=\{0,p,2p,\cdots,kp,kp+r,kp+p,\cdots\}$.
This implies that $\mathit{CFK}^\infty(T_{p,q})=\mathrm{St}((1,p-1)^k,1,r-1,1,p-r-1,\cdots)$.
Note that no entry in the staircase complex can exceed $p-1$.
Otherwise, such an entry could be assumed to appear in an even slot by the palindromicity and there would be $p$ consecutive numbers not in the semigroup, which is impossible since $p$ is a generator.
Then $[\![T_{p,q}]\!]=[(1,p-1)^k,(p-1,1)^k]+[1,r-1,1,p-r-1,\cdots]$ by Lemma~\ref{split}.
Because $r<p/2\Rightarrow r-1<p-r-1$, it follows from~\cite[Lemma 4.4]{1nn1} that $[1,r-1,1,p-r-1,\cdots]\gg[1,r-1,r-1,1]$.
\hfill$\Box$

\begin{Prop}\label{three torus knots} Suppose $p$ and $q$ are relatively prime positive integers with $4\leqslant p<q/2$ and $k$ is any positive integer.
Then $[\![T_{q,kq+p}\#-T_{p,q}\#-kT_{q,q+1}]\!]\gg[1,p-1,p-1,1]$ in $\mathcal{CFK}_\mathrm{alg}$.\end{Prop}

\textbf{Proof.}\begin{align*}&[\![T_{q,kq+p}\#-T_{p,q}\#-kT_{q,q+1}]\!]\\
=&[\![T_{q,kq+p}]\!]-[\![T_{p,q}]\!]-k[\![T_{q,q+1}]\!]\\
=&(k[1,q-1,q-1,1]+O_1)-(\lfloor q/p\rfloor[1,p-1,p-1,1]+O_2)-k([1,q-1,q-1,1]+O_3)\\
=&O_1-\lfloor q/p\rfloor[1,p-1,p-1,1]-O_2-kO_3.\end{align*}

Here $O_1\gg[1,p-1,p-1,1]$ by Lemma~\ref{half remainder}, $O_2$ is dominated by $[1,p-1,p-1,1]$ by Lemma~\ref{any remainder} and $O_3$ is dominated by $[1,p-1,p-1,1]$ by Lemma~\ref{remainder1}.

Finally the conclusion follows from the definition of $\gg$.
\hfill$\Box$

\emph{Remark.} This proposition implies that $a_1(K)=1$ and $a_2(K)\geqslant p-1$ for $K=T_{q,kq+p}\#-T_{p,q}\#-kT_{q,q+1}$ with $p,q,k$ as in the hypothesis.
It follows from~\cite[Proposition 3.12]{splitting} that there are knots with vanishing $\Upsilon$ invariant but arbitrarily large \emph{splitting concordance genus}~\cite[Definition 1.1]{splitting}.
In particular, there are knots with vanishing $\Upsilon$ invariant but arbitrarily large concordance genus, which was first obtained in~\cite{secondary0}.

Recall that~\cite{survey} defines the notion of \emph{stable equivalence} of complexes, which is finer than $\varepsilon$-equivalence.

\begin{Cor}Suppose $p$ and $q$ are relatively prime positive integers with $4\leqslant p<q/2$ and $k$ is any positive integer.
Then $\varepsilon(T_{q,kq+p}\#-T_{p,q}\#-kT_{q,q+1})=1$.
In particular, $\mathit{CFK}^\infty(T_{q,kq+p})$ is not stably equivalent to $\mathit{CFK}^\infty(T_{p,q}\#kT_{q,q+1})$.\end{Cor}

\emph{Remark.} The conclusion in the case $k=1$ is covered by~\cite[Theorem 1.2]{secondary2} (also see~\cite[Conjecture 5.3]{secondary1}).

In the proof of Proposition~\ref{three torus knots}, we know $O_1$ is dominated by $[1,q-1,q-1,1]$ by Lemma~\ref{any remainder}, even without the hypothesis $p<q/2$.
Meanwhile, $O_2$ and $O_3$ are dominated by $[1,p-1,p-1,1]$ and hence also dominated by $[1,q-1,q-1,1]$ by Lemma~\ref{a1a2greater}.
This immediately gives the following bound.

\begin{Lem}\label{upper bound}Suppose $p$ and $q$ are relatively prime positive integers with $4\leqslant p<q$ and $k$ is any positive integer.
Then $[\![T_{q,kq+p}\#-T_{p,q}\#-kT_{q,q+1}]\!]$ is dominated by $[1,q-1,q-1,1]$.\end{Lem}

\textbf{Proof of Theorem~\ref{main}.} Let $K_i$ be the knot $T_{q_i,k_iq_i+p_i}\#-T_{p_i,q_i}\#-k_iT_{q_i,q_i+1}$ for each positive integer $i$,
where $p_i$ and $q_i$ are relatively prime positive integers with $4\leqslant p_i<q_i/2$ and $k_i$ is any positive integer.
Any linear combination of the family $\{K_i\}_{i=1}^\infty$ obviously has vanishing $\Upsilon$ invariant by Theorem~\ref{expansion}.

We choose $p_i$ and $q_i$ such that $q_i\leqslant p_{i+1}$.
Then the family $\{K_i\}_{i=1}^\infty$ satisfies $[\![K_i]\!]\ll[\![K_{i+1}]\!]$ for each $i$,
because $[\![K_i]\!]$ is dominated by $[1,q_i-1,q_i-1,1]$ by Lemma~\ref{upper bound} and $[\![K_{i+1}]\!]\gg[1,p_i-1,p_i-1,1]$ by Proposition~\ref{three torus knots}.
Now $\{[\![K_i]\!]\}_{i=1}^\infty$ generates a subgroup isomorphic to $\mathbb{Z}^\infty$ in $\mathcal{CFK}_\mathrm{alg}$ by Lemma~\ref{domination implies independence}.
Therefore $\{K_i\}_{i=1}^\infty$ generates a subgroup isomorphic to $\mathbb{Z}^\infty$ in $\mathcal{C}$.
Any nontrivial linear combination of this family has nonvanishing $\varepsilon$ invariant, because it maps to a nontrivial linear combination of $\{[\![K_i]\!]\}_{i=1}^\infty$
under the homomorphism from $\mathcal{C}$ to $\mathcal{CFK}_\mathrm{alg}$ that maps any knot to its $\varepsilon$-equivalence class.
\hfill$\Box$

\begin{Ex}Take $p_i=3^i+1,q_i=2\cdot3^i+3,k_i=1$ for all $i$.
Then we have a family $\{T_{9,13}\#-T_{4,9}\#-T_{9,10},T_{21,31}\#-T_{10,21}\#-T_{21,22},T_{57,85}\#-T_{28,57}\#-T_{57,58},\cdots\}$,
which satisfies the conclusion of Theorem~\ref{main}.\end{Ex}

\end{document}